\documentclass[10pt,fleqn]{article}

\usepackage{porousmedia}
\usepackage{porousmedia-macros}

\title{Performant low-order matrix-free finite element kernels on GPU architectures}

\author[1]{Randolph~R.~Settgast
           \thanks{Corresponding author. \textit{E-mail address}: \texttt{settgast1@llnl.gov}}}
\author[2]{Yohann~Dudouit}
\author[1]{Nicola~Castelletto}
\author[3]{William~R.~Tobin}
\author[1]{Benjamin~C.~Corbett}
\author[4]{Sergey~Klevtsov}

\affil[1]{Atmospheric, Earth, and Energy Division, Lawrence Livermore National Laboratory, Livermore, CA, United States}
\affil[2]{Center for Applied Scientific Computing, Lawrence Livermore National Laboratory, Livermore, CA, United States}
\affil[3]{Applications, Simulations, and Quality Division, Lawrence Livermore National Laboratory, Livermore, CA, United States}
\affil[4]{Energy Resources Engineering, Stanford University, Stanford, United States}

\begin{document}

\maketitle

\begin{abstract}

Numerical methods such as the Finite Element Method (FEM) have been successfully adapted to utilize the computational power of GPU accelerators.
However, much of the effort around applying FEM to GPU’s has been focused on high-order FEM due to higher arithmetic intensity and order of accuracy.
For applications such as the simulation of subsurface processes, high levels of heterogeneity results in high-resolution grids characterized by highly discontinuous (cell-wise) material property fields.
Moreover, due to the significant uncertainties in the characterization of the domain of interest, e.g. geologic reservoirs, the benefits of high order accuracy are reduced, and low-order methods are typically employed.
In this study, we present a strategy for implementing highly performant low-order matrix-free FEM operator kernels in the context of the conjugate gradient (CG) method.
Performance results of matrix-free Laplace and isotropic elasticity operator kernels are presented and are shown to compare favorably to matrix-based SpMV operators on V100, A100, and MI250X GPUs.


\end{abstract}

\allowdisplaybreaks{}


\section{Introduction}
\label{sec:introduction}
The evolution of high performance computing architectures towards GPU accelerator based systems has shifted focus of many numerical method implementation towards data bandwidth considerations.
In many cases, numerical algorithms that were traditionally limited by FLOP/s count are now bandwidth bound.
The implementation of the finite element (FE) method for unstructured meshes is one such case.
The response of the FE community has been to move towards high order shape functions which achieve higher accuracy at a higher FLOP/s demand for a nearly fixed degree of freedom count and bandwidth demand.
This shift towards high order FE methods has been successful through efforts such as the Center for Efficient Exascale Discretizations (CEED) within the Department of Energy's (DOE) Exacale Computing Project (ECP). \cite{ceedbenchmarks,githubceedbenchmarks,ms32,ms34,ms38,ms39}

While a viable approach to many problems, the use of high order FE methods is not appropriate for all applications.
Simulations of many engineering and geoscience applications may involve material properties and state variables with high levels of uncertainty and heterogeneity.
For example in subsurface processes, this data is typically provided at a resolution representative of the heterogeneity/natural layering of the geologic formation as determined by a geophysical analysis.
The resulting fields of material properties are typically represented as a highly discontinuous set of values with an assumption of $C^0$ continuity within the tabular interval in which they are defined, and the use of higher-order discretizations does not practically offer the accuracy advantages over low-order approaches.
As such, low-order methods dominate a wide-range of simulation communities.

Application of a low order FE method typically involves the allocation and filling of a global sparse matrix $\Mat{A}$ of size ($ndo\!f \times nnzr$), where $ndo\!f$ is the global number of degrees of freedom in the problem, and $nnzr$ is the number of non-zeros in a row of the matrix.
For Cartesian grids, there are approximately ($27 \times ndo\!f\!sp$) number of non-zeros per row, where $ndo\!f\!sp$ is the number of degrees of freedom per support point.
For iterative Krylov methods, these large sparse matrices are used as a vector operator through a sparse matrix vector multiply.
In contrast, a matrix-free approach eliminates the need to store and manipulate the sparse global stiffness matrix, which results in a significant reduction in the amount of memory required to perform a simulation.
This is particularly important for hardware platforms in which the amount of available memory is the limiting factor on the number of degrees of freedom that may be solved in a simulation.
A second point of comparison between matrix-based and matrix-free operators is the speed of operator application.
In the case of a matrix-based operator, the speed of sparse matrix vector product has a low arithmetic intensity, i.e. is bound by the memory bandwidth (e.g. \cite{Zardoshti2017AdaptiveMultiplication}).
A matrix-free operator potentially allows a reduction in the bandwidth required to carry out the operator application, which can lead to faster runtime and better utilization of the GPU's capabilities in some cases.

In this paper we propose a novel strategy for implementing highly performant low-order matrix-free FE operator kernels.
Such kernels are used within a conjugate gradient (CG) iteration and are designed to take advantage of the computational power of GPU accelerators.
We performed a systematic performance analysis and comparison of the proposed matrix-free kernels against both matrix-based approaches and state-of-the-art performance results provided by the ECP CEED project \cite{ms32,ms34,ms38,ms39}.
Moreover, we highlight the key role played by the dot product implementation in the CG performance.
This analysis demonstrates that our low-order FEM implementation can achieve competitive performance on GPU architectures.

\section{Model problems}
\label{sec:model_problems}
We focus on the scalar Laplace's, vector Laplace's, and linear elasticity equations in three-dimensional domains.
The local matrices associated to an element $\Omega^e$ for such equations are defined as:
\begin{align}
  [\Mat{A}^e]_{ij}
  &
  = \int_{\Omega_e} \nabla \phi_i \cdot \nabla \phi_j \mathrm{d} \tensor{1}{x}
  && \text{(scalar Laplace)}, \\
  [\Mat{A}^e]_{3(i-1)+k,3(j-1)+l}
  &
  = \int_{\Omega_e}
    ( \tensor{1}{e}_k \otimes \nabla \phi_i )
    :
    ( \tensor{1}{e}_l \otimes \nabla \phi_j)
    \mathrm{d} \tensor{1}{x}
  && \text{(vector Laplace)}, \\
  [\Mat{A}^e]_{3(i-1)+k,3(j-1)+l}
  &
  = \int_{\Omega_e}
    ( \tensor{1}{e}_k \otimes \nabla \phi_i )
    :
    \tensor{2}{\sigma}( \tensor{1}{e}_l \otimes \nabla \phi_j)
    \mathrm{d} \tensor{1}{x}
  && \text{(linear elasticity)},
\end{align}
where
\begin{itemize}
  \item 
  $\{ \phi_i, \phi_j \}$ range over the scalar lowest-order nodal basis functions for cell $\Omega^e$, e.g. for an hexahedron $(i,j) \in \{1,2,\ldots,8\} \times \{ 1, 2, \ldots,8\}$;
  \item
  $\tensor{1}{e}_k$ are $\tensor{1}{e}_l$ denote standard-basis vectors in $\mathbb{R}^3$, with $k \in \{1,2,3 \}$ and $l \in \{1,2,3 \}$;
  \item
  $\tensor{2}{\sigma}$ is the second-order stress tensor that depends on the displacement gradient -- specifically on the symmetric gradient, i.e. strain tensor -- through an appropriate constitutive law;
  \item
  symbols `$\cdot$', `$:$', and `$\otimes$' denote the dot product, tensor contraction, and dyadic product, respectively, while $\nabla$ is the gradient operator. 
\end{itemize}
Evaluating integrals using an isoparametric transformation to a reference element and integrating numerically based on Gaussian quadrature, we obtain
\begin{align}
  [\Mat{A}^e]_{ij}
  &
  = \sum_q
    \left(
      \tensor{2}{J}_e^{-T} (\hat{\tensor{1}{x}}_q) \hat{\nabla} \hat{\phi}_i (\hat{\tensor{1}{x}}_q)
    \right)
    \cdot
    \left(
      \tensor{2}{J}_e^{-T} (\hat{\tensor{1}{x}}_q) \hat{\nabla} \hat{\phi}_j (\hat{\tensor{1}{x}}_q)
    \right)    
    w_q \determinant{\tensor{2}{J}_e (\hat{\tensor{1}{x}}_q)}
  && \text{(scalar Laplace)}, \label{eq:scalar_laplace_Ae} \\
  [\Mat{A}^e]_{3(i-1)+k,3(j-1)+l}
  &
  = \sum_q
    \left(
      \tensor{1}{e}_k
      \otimes
      \tensor{2}{J}_e^{-T} \hat{\nabla} \hat{\phi}_i (\hat{\tensor{1}{x}}_q)
    \right)
    :
    \left(
      \tensor{1}{e}_l
      \otimes
      \tensor{2}{J}_e^{-T} \hat{\nabla} \hat{\phi}_l (\hat{\tensor{1}{x}}_q)
    \right)    
    w_q \determinant{\tensor{2}{J}_e (\hat{\tensor{1}{x}}_q)}
  \nonumber \\
  &
  = \sum_q
    \tensor{1}{e}_k
    \cdot
    \left(
      \tensor{1}{e}_l
      \otimes
      \hat{\nabla} \hat{\phi}_j (\hat{\tensor{1}{x}}_q)
    \right)    
    \tensor{2}{J}_e^{-1}
    \tensor{2}{J}_e^{-T}     
    \hat{\nabla} \hat{\phi}_i (\hat{\tensor{1}{x}}_q)
    w_q \determinant{\tensor{2}{J}_e (\hat{\tensor{1}{x}}_q)}
  && \text{(vector Laplace)}, \label{eq:vector_laplace_Ae} \\
  [\Mat{A}^e]_{3(i-1)+k,3(j-1)+l}
  &
  = \sum_q
    \left(
      \tensor{1}{e}_k
      \otimes
      \tensor{2}{J}_e^{-T} \hat{\nabla} \hat{\phi}_i (\hat{\tensor{1}{x}}_q)
    \right)
    :
    \tensor{2}{\sigma}
    \left(
      \tensor{1}{e}_l
      \otimes
      \tensor{2}{J}_e^{-T} \hat{\nabla} \hat{\phi}_l (\hat{\tensor{1}{x}}_q)
    \right)    
    w_q \determinant{\tensor{2}{J}_e (\hat{\tensor{1}{x}}_q)}
  \nonumber \\
  &
  = \sum_q
    \tensor{1}{e}_k
    \cdot
    \tensor{2}{\sigma}
    \left(
      \tensor{1}{e}_l
      \otimes
      \hat{\nabla} \hat{\phi}_j (\hat{\tensor{1}{x}}_q)
      \tensor{2}{J}_e^{-1}
    \right)    
    \tensor{2}{J}_e^{-T}     
    \hat{\nabla} \hat{\phi}_i (\hat{\tensor{1}{x}}_q)
    w_q \determinant{\tensor{2}{J}_e (\hat{\tensor{1}{x}}_q)}
  && \text{(linear elasticity)}, \label{eq:elasticity_Ae}
\end{align}
where
\begin{itemize}
  \item 
  $\hat{\tensor{1}{x}}_q$ and $w_q$ denote quadrature point coordinates in the reference element and weights, respectively;
  \item
  $\tensor{2}{J}_e$ is the Jacobian transformation matrix from the reference element to element $\Omega^e$, and $\determinant{\tensor{2}{J}_e}$ its determinant;
  \item
  $\hat{\nabla}$ denotes the reference space gradient.
\end{itemize}
For notational convenience, we will not explicitly indicate from now on that terms in the reference element are evaluated at $\hat{\tensor{1}{x}}_q$.

\section{Matrix-free finite element operator}
\label{sec:matrix_free_FE_operator}
Solving a linear system $A \Vec{x} = \Vec{b}$ by a Krylov subspace method requires the application of the operator $\Mat{A}$ to one (or more) vectors at each iteration.
In matrix based methods, the application of matrix $\Mat{A}$ to a vector is accomplished by storing the operator in a sparse matrix format, and performing a sparse matrix vector multiplication.
%
In a matrix-free method, the sparse matrix is not formed, and the operator is applied directly to the input vector in a manner similar to that which is used to construct $\Mat{A}$.
Indeed, applying matrix $\Mat{A}$ to a vector $\Vec{u}$ to compute vector $\Vec{v}$ in a matrix-free fashion requires merging the matrix assembly into the multiplication and consists of three steps for each element $\Omega^e$:
\begin{enumerate}
  \item
  Reading the local components $\Vec{u}^e$ from the global input vector $\Vec{u}$;
  \item
  Evaluating the local matrix multiplication $\Vec{v}^e = \Mat{A}^e \Vec{u}^e$;
  \item
  Assembling the local components of $\Vec{v}^e$ into the global output vector $\Vec{v}$.
\end{enumerate}
Using expressions \eqref{eq:scalar_laplace_Ae}-\eqref{eq:elasticity_Ae} into the local multiplication yields
\begin{align}
  [\mathbf{v}^e]_{i}
  &
  = \sum_q \sum_j
    \left(
      \tensor{2}{J}_e^{-T} \hat{\nabla} \hat{\phi}_i
    \right)
    \cdot
    \left(
      \tensor{2}{J}_e^{-T} \hat{\nabla} \hat{\phi}_j
    \right)    
    w_q \determinant{\tensor{2}{J}_e}
    [\mathbf{u}^e]_{j}
  \nonumber \\
  &
  = \sum_q
    \hat{\nabla} \hat{\phi}_i
    \cdot
    \tensor{2}{J}_e^{-1}
    \left(
      \tensor{2}{J}_e^{-T}
      \sum_j
        \hat{\nabla} \hat{\phi}_j [\mathbf{u}^e]_{j}
    \right)    
    w_q \determinant{\tensor{2}{J}_e}
  && \text{(scalar Laplace)}, \label{eq:scalar_laplace_Ae_ue} \\
  [\mathbf{v}^e]_{3(i-1)+k}
  &
  = \sum_q \sum_j \sum_l
    \tensor{1}{e}_k
    \cdot
    \left(
      \tensor{1}{e}_l
      \otimes
      \hat{\nabla} \hat{\phi}_j
    \right)    
    \tensor{2}{J}^{-1}
    \tensor{2}{J}^{-T}     
    \hat{\nabla} \hat{\phi}_i
    w_q \determinant{\tensor{2}{J}_e}
    [\mathbf{u}^e]_{3(j-1)+l}
  \nonumber \\
  &
  = \sum_q
    \tensor{1}{e}_k
    \cdot
    \left(
      \sum_j
      \tensor{1}{u}^e_j
      \otimes
      \hat{\nabla} \hat{\phi}_j
    \right)    
    \tensor{2}{J}^{-1}
    \tensor{2}{J}^{-T}     
    \hat{\nabla} \hat{\phi}_i\
    w_q \determinant{\tensor{2}{J}_e}
  && \text{(vector Laplace)},  \label{eq:vector_laplace_Ae_ue} \\
  [\mathbf{v}^e]_{3(i-1)+k}
  &
  = \sum_q \sum_j \sum_l
    \tensor{1}{e}_k
    \cdot
    \tensor{2}{\sigma}
    \left(
      \tensor{1}{e}_l
      \otimes
      \hat{\nabla} \hat{\phi}_j
      \tensor{2}{J}^{-1}
    \right)    
    \tensor{2}{J}^{-T}     
    \hat{\nabla} \hat{\phi}_i
    w_q \determinant{\tensor{2}{J}_e}
    [\mathbf{u}^e]_{3(j-1)+l}
  \nonumber \\
  &
  = \sum_q
    \tensor{1}{e}_k
    \cdot
    \tensor{2}{\sigma}
    \left(
      \left(
        \sum_j
        \tensor{1}{u}^e_j
        \otimes
        \hat{\nabla} \hat{\phi}_j   
      \right)
      \tensor{2}{J}^{-1}
    \right) 
    \tensor{2}{J}^{-T}     
    \hat{\nabla} \hat{\phi}_i\
    w_q \determinant{\tensor{2}{J}_e}
  && \text{(linear elasticity)}.  \label{eq:elasticity_Ae_ue}
\end{align}
Note in Eqs.~\eqref{eq:vector_laplace_Ae_ue} and \eqref{eq:elasticity_Ae_ue} that vector $\tensor{1}{u}^e_j = \sum_l \tensor{1}{e}_l [\mathbf{u}^e]_{3(j-1)+l}$ contains the three components of the algebraic vector $\Vec{u}^e$ associated with node $j$.

While matrix-based methods have traditionally been the standard for low-order finite element problems, matrix-free methods have characteristics that make them an interesting option.
The most significant advantage of matrix-free methods is eliminating the need to store and fill the sparse global matrix.
In a FE method the sparse global matrix typically represents a significant percentage of total memory required to run a simulation and removal of the allocation requirement allows a much larger problem to be run given a fixed total memory constraint.

In addition to the reduced memory overhead of a matrix-free method, in some cases the amount of data transfer costs required by a matrix-free kernel are less than the associated cost to transfer the sparse matrix into the kernel.
In such cases, it is possible that the matrix-free method may outperform a matrix-based method in terms of execution speed.

One more advantage that is highlighted here is the removal of the requirement to fill the sparse matrix. 
When implementing a matrix-based approach, the algorithm typically involves a loop over elements where a dense element local stiffness matrix is computed using some quadrature rule, see Eqs.~\eqref{eq:scalar_laplace_Ae}-\eqref{eq:elasticity_Ae}.
The values of this matrix are then added to their appropriate location in the sparse global matrix.
This process of adding local matrix contributions to the global sparse matrix typically involves some sort of search for the global column/dof location in the sparse matrix followed by an atomic add operation.
In linear problems the cost of the fill operation may be insignificant to the cost of solving the linear system as it is performed only once.
However for a nonlinear problem, the fill operation will occur once per Newton iteration, and the sparse matrix fill may account for a significant portion of the overall timestep.

In this work, we will present results from a simple matrix-free element operator as described in Algorithm~\ref{alg:local_stiffness_matrix_free_randy} for the elastotastics case.
We note that the approach used by the ECP CEED project consists in precomputing and storing the determinant of the Jacobian transformation and its inverse all the values at quadrature points.
We will refer in the rest of this article to this matrix-free approach as "partial assembly".

\begin{algorithm}
\caption{Calculate action of stiffness operator (matrix-free) for quasi-static mechanics}
\label{alg:local_stiffness_matrix_free_randy}
\renewcommand{\algorithmicrequire}{\textbf{Input:}}
\renewcommand{\algorithmicensure}{\textbf{Output:}}
\begin{algorithmic}[1]
\Require 
    \Statex $\Vec{x}$: global support point coordinates
    \Statex $\Vec{u}$: global input vector
\Ensure  
    \Statex $\Vec{v}$: global output vector 
    \newline 
\State $\Vec{x}^e = \{ \tensor{1}{x}_j^e \}_{j = 1}^{\text{\# nodes}} \leftarrow \Vec{x}$,
       $\Vec{u}^e = \{ \tensor{1}{u}_j^e \}_{j = 1}^{\text{\# nodes}}\leftarrow \Vec{u}$
       \Comment{Gather/Load data into element local storage }
\State $\Vec{v}^e = \{ \tensor{1}{v}_j^e \}_{j = 1}^{\text{\# nodes}}$
       \Comment{Allocate local storage for output vector}
\For{each quadrature point $q$}
  
  \State $ \tensor{2}{J}_e
         = \sum_j^{\text{\# nodes}}
           \tensor{1}{x}_j^e \otimes \hat{\nabla} \hat{\phi}_j $
         \Comment{Calculate Jacobian transformation}
         
  \State $\tensor{2}{J}_e^{-1} = \text{inverse}\left( \tensor{2}{J}_e \right)$,
         $\determinant{\tensor{2}{J}_e} = \text{determinant}\left(\tensor{2}{J}_e \right) $
         \Comment{Calculate inverse Jacobian transformation and determinant}

  \State $ \nabla \tensor{1}{u}
         = \left(
           \sum_j^{\text{\# nodes}}
             \tensor{1}{u}_j^e
             \otimes
             \hat{\nabla} \hat{\phi}_j
           \right)
           \tensor{2}{J}_e^{-1}$
         \Comment{Calculate gradient of input variable}
  \State $\tensor{2}{\sigma} = \tensor{2}{\sigma}(\nabla \tensor{1}{u})$ 
         \Comment{Calculate stress tensor}
  \State $\tensor{2}{\sigma} \timeseq w_q \determinant{\tensor{2}{J}_e}$ 
         \Comment{Scale stress by Jacobian determinant times quadrature weight}      
  \State $\tensor{1}{P} =  \tensor{2}{\sigma} \tensor{2}{J}_e^{-T}$
         \Comment{Multiply stress by test function gradients (step 1)}
  \State $\tensor{1}{v}_i^e \pluseq \tensor{2}{P} \hat{\nabla} \hat{\phi}_i$, with $i \in \{1, \ldots, \text{\# nodes} \}$
         \Comment{Multiply stress by test function gradients (step 2)}
\EndFor
\State $\Vec{v}^e \rightarrow \Vec{v}$
       \Comment{Add local contributions to global output vector} 
\end{algorithmic}
\end{algorithm}

\section{Theoretical kernel memory transfer and idealized performance}
\label{sec:theoretical_kernel_analysis}
As a basis of comparison between operator kernels, we propose theoretical memory requirements and throughput for the cases of scalar Laplace and isotropic linear elasticity.
The operators for sparse matrix-vector product (SpMV), matrix-free, and partial assembly are evaluated.
The memory and throughputs are intended to place a theoretical bound on the maximum theoretical performance for each operator applied to scalar Laplace and isotropic linear elastic problems.
The method by which the quantities are calculated is to take a problem of one million elements, and calculate the data arrays that are loaded in order to execute the operator kernel.
Then the published bandwidth of each architecture is used to calculate a "memory transfer time", which is also referred to as "Speed of Light" in the community.
Finally the number of degrees of freedom is divided by the "memory transfer time" to estimate a theoretical throughput.

\begin{table}[htbp]
  \small
  \centering
  \caption{ Theoretical throughput for SpMV. 
  \textmd{We present theoretical memory transfer requirements executing a sparse matrix vector product (SpMV) kernel for Scalar Laplace, and mechanics.
  Additionally, memory transfer times and throughput are presented assuming published bandwidths for the NVIDIA V100, NVIDIA A100, and AMD MI250X.} }
    \begin{tabular}{l c c}
    \toprule
    Quantity                & Scalar Laplace   & Mechanics \\ 
    \midrule
    number of elements          & \multicolumn{2}{c}{1,000,000}\\
    number of nodes             & \multicolumn{2}{c}{1,030,301}\\
    memory / non-zero (bytes)   & \multicolumn{2}{c}{12}\\
    \midrule
    num dof/support point  & 1          & 3 \\
    rows (num dof)         & 1,030,301  & 3,090,903 \\
    nnz/row                & 27         & 81 \\
    \midrule
    matrix memory (MB)     & 343        & 3,004 \\
    vector memory (MB)     & 16         & 49 \\
    total memory (MB)      & 359        & 3,079 \\ 
    \midrule
    V100 transfer time (ms)   & 0.40	& 3.4 \\
    A100 transfer time (ms)   & 0.19	& 1.6 \\
    MI250X transfer time (ms) & 0.22	& 1.9 \\
    \midrule
    V100 throughput (GDof/s)   & 2.6	& 0.90 \\
    A100 throughput (GDof/s)   & 5.6	& 1.9 \\
    MI250X throughput (GDof/s) & 4.7	& 1.6 \\
    \bottomrule
    \end{tabular}%
  \label{tab:SpMVMult_memoryRequirements}%
\end{table}%

We begin by providing these data for a sparse matrix vector product (see Table \ref{tab:SpMVMult_memoryRequirements}) as a reference.
It is important to note that linear elasticity contains three times the number of rows when compared to scalar Laplace.
In term of non-zeros, scalar Laplace contains 27 non-zeros per row, whereas linear elasticity contains 81 non-zeros per row.
This leads to a significantly larger memory requirements for linear elasticity when comparing problems with equivalent number of elements.



For matrix-free algorithms, we consider two idealized scenarios for memory transfers, one with a perfect memory caching, and one with no-caching at all.
The perfect caching scenario assumes that the cache infrastructure will hide the cost to read multiple times the same data.
The no-caching scenario assumes that each data that needs to be read multiple times will not benefit from any caching.
These two scenarios give us an idealized performance range for each finite element operator and for each GPU architecture.
In Table~\ref{tab:matrixFreeMemoryRequirements} we present these idealized performance for the matrix-free algorithm that we are proposing.
In Table~\ref{tab:partialAssemblyMemoryRequirements} we present these idealized performance for the "partial assembly" algorithm used by the ECP CEED project.

\begin{table}
  \small
  \centering
  \caption{ Theoretical throughput for matrix-free FE operator. 
  \textmd{We present theoretical memory transfer requirements executing a matrix-free FE operator kernel for Scalar Laplace, and mechanics assuming isotropic elasticity.
  Nodal memory ranges are intended to reflect different assumptions of read efficiency; from each nodal variable being read once, to each nodal variable being read once for each element it is connected to.
  Additionally, memory transfer times and throughput are presented assuming published bandwidths for the NVIDIA V100, NVIDIA A100, and AMD MI250X.} }
  \begin{tabular}{ l c c }
    \toprule
            & Laplace & isotropic elasticity \\
    \midrule
    node map (MB)                           & 32   & 32   \\
    cell constant \textbf{C} (MB)           &  0   & 16   \\
    \# data values / quadrature pt          &  0   &  0   \\
    quadrature storage (MB)                 & 0    &  0   \\
    \midrule
    read nodal position (MB)             & 25 - 192  & 25 - 192 \\
    input/output vectors (MB)            & 25 - 192  & 74 - 576  \\
    \midrule
    total memory (MB)                    & 81 - 416 & 147 - 816  \\
    \midrule
    V100 transfer time (ms)   & 0.090 - 0.46	& 0.16 - 0.91 \\
    A100 transfer time (ms)   & 0.042 - 0.22	& 0.076 - 0.42 \\
    MI250X transfer time (ms) & 0.049 - 0.25	& 0.090 - 0.50 \\
    \midrule
    V100 throughput (GDof/s)   & 11 - 2.2	& 19 - 3.4 \\
    A100 throughput (GDof/s)   & 25 - 4.8	& 41 - 7.3 \\
    MI250X throughput (GDof/s) & 21 - 4.1	& 34 - 6.2 \\
    \bottomrule
  \end{tabular}%
  \label{tab:matrixFreeMemoryRequirements}%
\end{table}%

\begin{table}
  \centering
  \small
  \caption{ Theoretical throughput for FE operator using "partial assembly". 
  \textmd{We present theoretical memory transfer requirements executing a partial assembly FE operator kernel (need reference) for Scalar Laplace, and mechanics assuming isotropic elasticity.
  Nodal memory ranges are intended to reflect different assumptions of read efficiency; from each nodal variable being read once, to each nodal variable being read once for each element it is connected to.
  Additionally, memory transfer times and throughput are presented assuming published bandwidths for the NVIDIA V100, NVIDIA A100, and AMD MI250X.} }
  \begin{tabular}{ l c c }
    \toprule
            & Laplace & isotropic elasticity \\ 
    \midrule
    node map (MB)                           & 32   &    32   \\
    cell constant \textbf{C} (MB)           &  0   &    16   \\
    \# data values / quadrature pt          &  6   &    21   \\
    quadrature storage (MB)                 & 384  & 1,344   \\
    \midrule
    input/output vectors (MB)               & 25 - 192  & 74 - 576   \\
    \midrule
    total memory (MB)                       & 441 - 608&	1,450 - 1,952  \\
    \midrule
    V100 transfer time (ms)   & 0.49 - 0.68 & 1.6 - 2.2 \\
    A100 transfer time (ms)   & 0.23 - 0.31 & 0.75 - 1.0 \\
    MI250X transfer time (ms) & 0.27 - 0.37 & 0.89 - 1.2 \\
    \midrule
    V100 throughput (GDof/s)   & 2.1 - 1.5 & 1.9 - 1.4 \\
    A100 throughput (GDof/s)   & 4.5 - 3.3 & 4.1 - 3.1 \\
    MI250X throughput (GDof/s) & 3.8 - 2.8 & 3.5 - 2.6 \\
    \bottomrule
  \end{tabular}%
  \label{tab:partialAssemblyMemoryRequirements}%
\end{table}%

Comparing idealized performance of a sparse matrix-vector product with our matrix-free algorithm, we observe that our algorithm is more favorable to vector problems than scalar problems.
Indeed, the amount of data per degree-of-freedom is generally lower for vector problems than for scalar problems due to increased data reuse in vector problems.
Overall, we observe that the idealized sparse matrix-vector performance is lower than the worst performance scenario for our matrix-free algorithm (see Table~\ref{tab:SpMVMult_memoryRequirements} and Table~\ref{tab:matrixFreeMemoryRequirements} ).

Comparing idealized performance of our matrix-free algorithm with the "partial assembly" algorithm, we observe that the worst case performance scenario for our algorithm is always above the best case performance scenario for "partial assembly" (see Table~\ref{tab:matrixFreeMemoryRequirements} and Table~\ref{tab:partialAssemblyMemoryRequirements}).
In the best case scenario of perfect caching, we observe that our matrix-free algorithm is a full order of magnitude faster than "partial assembly" for isotropic elasticity.

\begin{table}
  \centering
  \small
  \caption{Theoretical memory transfer size (MB) for Conjugate Gradient kernels (1e6 elements) }
    \begin{tabular}{ l c c c }
    \toprule
    Operation                                                   & Laplace   & \multicolumn{2}{c}{Mechanics} \\ 
    \midrule
    ddot $\Vec{p}_k^T \Mat{A} \Vec{p}_k$                             & 16.5    & \multicolumn{2}{c}{49.5} \\
    \midrule
    axpy $\Vec{x}_{k+1}=\Vec{x}_k+\alpha_k \Vec{p}_k$       & 16.5    & \multicolumn{2}{c}{49.5} \\
    \midrule
    axpy $\Vec{r}_{k+1}=\Vec{r}_k-\alpha_k \Mat{A}\Vec{p}_k$      & \multirow{2}{*}{ 16.5 }  & \multicolumn{2}{c}{\multirow{2}{*}{ 49.5 }} \\
    ddot $\Vec{r}_{k+1}^T \Vec{r}_{k+1}$                        \\
    \midrule
    axpy $\Vec{p}_{k+1}=\Vec{r}_{k+1}+\beta_k \Vec{p}_k$        & 16.5    & \multicolumn{2}{c}{49.5} \\ 
    \bottomrule
    \end{tabular}%
  \label{tab:cgMemoryRequirements}%
\end{table}%

\section{Code design and efficient GPU implementation}
\label{sec:code_design_GPU_implementation}
In this section, we discuss the critical aspects of code design that contribute to the observed high performance in matrix-free low-order FEM operator kernels.
The key elements include a generic programming approach, the use of \verb|static constexpr| methods and data-free classes for efficient data management and GPU utilization, and a moderate use of meta-programming techniques for loop unrolling and forcing constant inlining.

    \subsection{Generic programming approach}
    A generic programming approach has been employed in the code design to enable flexibility and reusability.
    By using templates and abstract data types, a wide range of problem types and computational kernels are accommodated without making changes to the core code.
    This approach promotes code maintainability and allows for easier adaptation to new algorithms or hardware architectures.

    Specifically in this work, FE operator kernel classes are templated on the finite element class, and the constitutive model.
    The compile time knowledge of the finite element class provides the compiler with the number of support points, the quadrature rule, and device callable inline functions for common finite element operations that are specialized for the finite element class.
    Compile time knowledge of the constitutive model type provides device callable constitutive stiffness operators.

    Multidimensional array types from the LvArray library 
    include in their template specification the underlying storage permutation, which allows optimization of memory access patterns to facilitate performant operations on GPU without any platform-specific kernel modification.
    
    Optimizing memory access patterns in this way has proven critical for achieving high performance on the GPU, as it enables efficient use of memory resources and reduces the potential for memory access bottlenecks.

    \subsection{\texttt{static constexpr} methods and data-free classes}
    The code relies heavily on static \verb|constexpr| methods, which enable computation at compile time, resulting in highly optimized code.
    The use of data-free classes containing only \verb|static constexpr| methods allows us to store information within the type itself rather than in separate data structures.
    This design choice reduces memory requirements and simplifies data management.

    Moreover, data-free classes facilitate efficient GPU utilization by enabling the transfer of classes containing only \verb|static constexpr| methods to the GPU for free, without any data movement.
    Since all the information is contained within the type itself, GPU resources can be allocated more effectively, and performance is further enhanced.

    \subsection{Meta-programming for loop unrolling and inlining}
    Meta-programming techniques are utilized to ensure efficient execution of our FEM operator kernels on GPUs.
    By employing compile-time loop unrolling and inlining, we eliminate the overhead associated with traditional loop constructs, which can significantly impact the performance of GPU-based computations.
    This optimization allows the compilers to use registers instead of local/scratch memory for arrays for instance.
    
    Such meta-programming was critical on AMD hardware where a comparison with \verb|#pragma unroll| showed significant performance improvement.
    Similarly, on NVIDIA hardware our meta-programming approach and leveraging \verb|#pragma unroll| resulted in similar performances.

    Aggressive compile-time recursive template instantiation and inlining of multidimensional array accessors results in these operations compiling down to direct memory accesses, allowing device compilers to fully optimize memory accesses with no uncertainty.

\section{Numerical results}
\label{sec:numerical_results}
In this section, we present performance results for matrix-free low-order FEM operator kernels that employ the proposed guidelines for efficient GPU implementations.
The results for the scalar Laplace and isotropic elasticity problems are shown for Nvidia V100, Nvidia A100, and AMD MI250X GPU's.

The CEED benchmark suite \cite{ceedbenchmarks, githubceedbenchmarks} is composed of a range of benchmark problems (BPs) that target different aspects of computational kernels commonly found in high-order finite element applications.
These BPs are designed to isolate specific operations, allowing researchers and developers to measure and analyze the performance of individual components of their implementations.

\begin{figure}[h!]
  \small
  \centering
  \if \generateTikzFigures 1
    \input{plotData/scalar_laplace_throughput}
  \else
    \includegraphics[scale=1]{./main-figure\theFigureCounter.pdf}
    \stepcounter{FigureCounter}
  \fi 
  \caption{
    Scalar Laplace - Conjugate gradient throughput.
    \textmd{We observe higher performance on AMD MI250X hardware than on NVIDIA V100 hardware.
    The NVIDIA V100 runs out of memory before achieving peak performance, on the other hand peak performance is achieved on AMD MI250X between 10 to 50 million elements.
    Strong scalability is slightly better on AMD MI250X than NVIDIA V100.
    The NVIDIA A100 outperforms the AMD MI250X on the BP5 benchmark despite having lower floating point operations potential and a similar memory bandwidth.
    This matrix-free algorithm achieves state-of-the-art performance for all these GPUs on the ECP CEED Bake-Off problem: BP5 \cite{ms32,ms34,ms38,ms39}.}
  }
  \label{fig:ScalarLaplaceThroughput}
\end{figure}

\begin{figure}[h!]
  \small
  \centering
  \if \generateTikzFigures 1
    \input{plotData/scalar_laplace_metrics_MI250X}
  \else
    \includegraphics[scale=1]{./main-figure\theFigureCounter.pdf}
    \stepcounter{FigureCounter}
  \fi 
  \caption{
    Scalar Laplace Kernel - MI250X performance metrics.
    \textmd{We observe high cache hit rates, and the plateauing of the global device memory bandwidth and FLOPs at approximately 40\% of their peak performance.}
  }
  \label{fig:ScalarLaplaceMetrics}
\end{figure}

\begin{figure}[h!]
  \small
  \centering
  \if \generateTikzFigures 1
    \input{plotData/scalar_laplace_kernel_throughput}
  \else
    \includegraphics[scale=1]{./main-figure\theFigureCounter.pdf}
    \stepcounter{FigureCounter}
  \fi 
  \caption{
    Scalar Laplace - Kernel throughput.
    \textmd{
        We observe that the NVIDIA V100 and AMD MI250X achieve performance beyond the maximum achievable performance for SpMV on their respective architectures.}
  }
  \label{fig:ScalarLaplaceKernelThroughput}
\end{figure}

\begin{figure}[h!]
  \small
  \centering
  \if \generateTikzFigures 1
    \input{plotData/vector_laplace_throughput}
  \else
    \includegraphics[scale=1]{./main-figure\theFigureCounter.pdf}
    \stepcounter{FigureCounter}
  \fi 
  \caption{
    Vector Laplace - Conjugate gradient throughput.
    \textmd{We observe similar but higher performance behavior for Vector Laplace than Scalar Laplace.
    This result is consistent with the lower amount of data moved per degree of freedom for the Vector Laplace operator.
    This algorithm achieves state-of-the-art performance on the ECP CEED Back-Off problem: BP6 \cite{ms32,ms34,ms38,ms39}.}
  }
  \label{fig:VectorLaplaceThroughput}
\end{figure}

\begin{figure}[h!]
  \small
  \centering
  \if \generateTikzFigures 1
    \input{plotData/solid_mech_throughput}
  \else
    \includegraphics[scale=1]{./main-figure\theFigureCounter.pdf}
    \stepcounter{FigureCounter}
  \fi 
  \caption{
    Solid mechanics - Conjugate gradient throughput.
    \textmd{The performance behavior for solid mechanics is very similar to the performance of the vector Laplace operator despite having more computation and slightly more data movement.
    The NVIDIA A100 significantly outperforms the AMD MI250X GPU, however peak performance is only achieved between 100 and 200 million degrees-of-freedom on NVIDIA A100.}
  }
  \label{fig:SolidMechThroughput}
\end{figure}

In this study, we specifically chose to focus on the ECP CEED benchmark problems BP5 (see Figure~\ref{fig:ScalarLaplaceThroughput}) and BP6 (see Figure~\ref{fig:VectorLaplaceThroughput}) for first order finite element method, as they are the most relevant to the applications we are addressing.
BP5 and BP6 correspond to the scalar and vectorial diffusion problems, respectively, and both involve the use of quadrature rules that are commonly employed in the types of applications we are targeting.
We also present results for a linear elasticity operator (see Figure~\ref{fig:SolidMechThroughput}) as a target benchmark problem using first order finite elements.

By selecting BP5 and BP6 as benchmark problems, we can assess the performance of matrix-free low-order FEM operator kernels in the context of the most widely used quadrature rules for our specific applications.
This choice enables us to draw more meaningful conclusions about the effectiveness and efficiency of our implementation, ensuring that the results are directly applicable to real-world problems.
By comparing our low-order finite element kernels with those of the ECP CEED project, we can identify the strengths and weaknesses of our implementation.

Low-order matrix-free finite element methods rely more on cache efficiency than high-order matrix-free methods due to the differences in their computational characteristics and memory access patterns.
This is particularly critical on GPU architectures, where memory bandwidth and efficient use of cache are key factors in achieving high performance.
However, their lower memory requirements also offer certain advantages that can help optimize cache usage and memory transactions.

\begin{itemize}
    \item Arithmetic intensity: Low-order finite element methods have lower arithmetic intensity compared to high-order methods. In low-order methods, fewer floating-point operations are performed per degree of freedom (DOF) than in high-order methods. On GPUs, which are designed to handle high arithmetic intensity, the low arithmetic intensity of low-order methods can result in underutilization of GPU resources, making efficient cache usage crucial to achieving high performance.

    \item Data locality: High-order methods have better data locality than low-order methods because the computations involve more DOFs within an element. In contrast, low-order methods involve fewer DOFs per element, leading to a lower data reuse rate and a higher reliance on cache efficiency to maintain performance on GPUs.

    \item Memory access patterns: In low-order methods, the memory access patterns tend to be more irregular than in high-order methods. Irregular memory access patterns can lead to cache misses and less efficient cache utilization on GPUs, where coalesced memory access is critical for optimal performance.

    \item GPU-specific challenges: Efficiently utilizing the GPU's memory hierarchy is essential for achieving high performance in low-order matrix-free finite element implementations. This involves optimizing shared memory usage, minimizing global memory accesses, and maximizing occupancy to hide memory latency. These challenges are more pronounced for low-order methods due to their lower arithmetic intensity and reduced data locality.

    \item Single thread per finite element: Low-order methods have lower memory requirements, which allows using a single thread per finite element. This approach eliminates the need for shared or local memory and relies solely on registers for storing the required data during computation.
    Registers are the fastest memory hierarchy available on GPUs, and utilizing them efficiently can significantly boost the performance of low-order finite element methods. 
    By using a single thread per finite element and leveraging register memory, memory transactions can be maximized, enabling the cache infrastructure to work optimally and achieving peak performance.
\end{itemize}

In summary, low-order matrix-free finite element methods rely more on cache efficiency than high-order methods due to their lower arithmetic intensity, reduced data locality, and irregular memory access patterns (see Figure~\ref{fig:ScalarLaplaceMetrics}).
This is particularly critical on GPU architectures, where memory bandwidth and efficient use of cache are key factors in achieving high performance.
However, the lower memory requirements of low-order methods also offer certain advantages, such as the ability to use a single thread per finite element and leverage register memory, which can help optimize cache usage and memory transactions.

\begin{figure}[h!]
  \small
  \centering
  \if \generateTikzFigures 1
    \input{plotData/CGvsOperator}
  \else
    \includegraphics[scale=1]{./main-figure\theFigureCounter.pdf}
    \stepcounter{FigureCounter}
  \fi 
  \caption{
    Throughput comparison of a Conjugate Gradient iteration with the sole operator kernel  - AMD MI250X.
    \textmd{We observe a significant performance discrepancy between the performance of the solver against the performance of the finite element operator alone, impacting significantly strong scalability and peak performance.}
  }
  \label{fig:CGvsOperator}
\end{figure}

\begin{figure}[h!]
  \small
  \centering
  \if \generateTikzFigures 1
    \input{plotData/CG_dot}
  \else
    \includegraphics[scale=1]{./main-figure\theFigureCounter.pdf}
    \stepcounter{FigureCounter}
  \fi 
  \caption{
    Conjugate gradient sensitivity to "dot" product implementation on AMD MI250X.
    \textmd{The "dot" product may seem like an harmless function, but its implementation has dramatic impact both on peak performance and strong scalability.}
  }
  \label{fig:Dot}
\end{figure}

\begin{figure}[h!]
  \small
  \centering
  \if \generateTikzFigures 1
    \input{plotData/solid_mech_kernel_throughput}
  \else
    \includegraphics[scale=1]{./main-figure\theFigureCounter.pdf}
    \stepcounter{FigureCounter}
  \fi 
  \caption{
    Solid mechanics - Kernel throughput.
    \textmd{The performance behavior on the sole solid mechanics kernel tells us a different story than a whole conjugate gradient iteration.
    Despite having similar peak bandwidth, the NVIDIA A100 and AMD MI250X have very different performance behavior.
    We also observe that contrary to the whole conjugate gradient iteration, kernel operator peak performance is achieved on all GPU architectures starting from one million elements.}
  }
  \label{fig:SolidMechKernelThroughput}
\end{figure}

The difference in throughput measurements between the sole finite element operator and the whole CG iteration (see Figure~\ref{fig:CGvsOperator}) is particularly noticeable when comparing the AMD MI250X and the NVIDIA A100 GPUs (see Figure~\ref{fig:SolidMechKernelThroughput} against Figure~\ref{fig:SolidMechThroughput}).
The performance discrepancy can be attributed to the reduction implementation in the dot product operation, which has a significant impact on the overall CG iteration performance (see Figure~\ref{fig:Dot})

\begin{figure}
    \centering
    \includegraphics[width=0.48\textwidth]{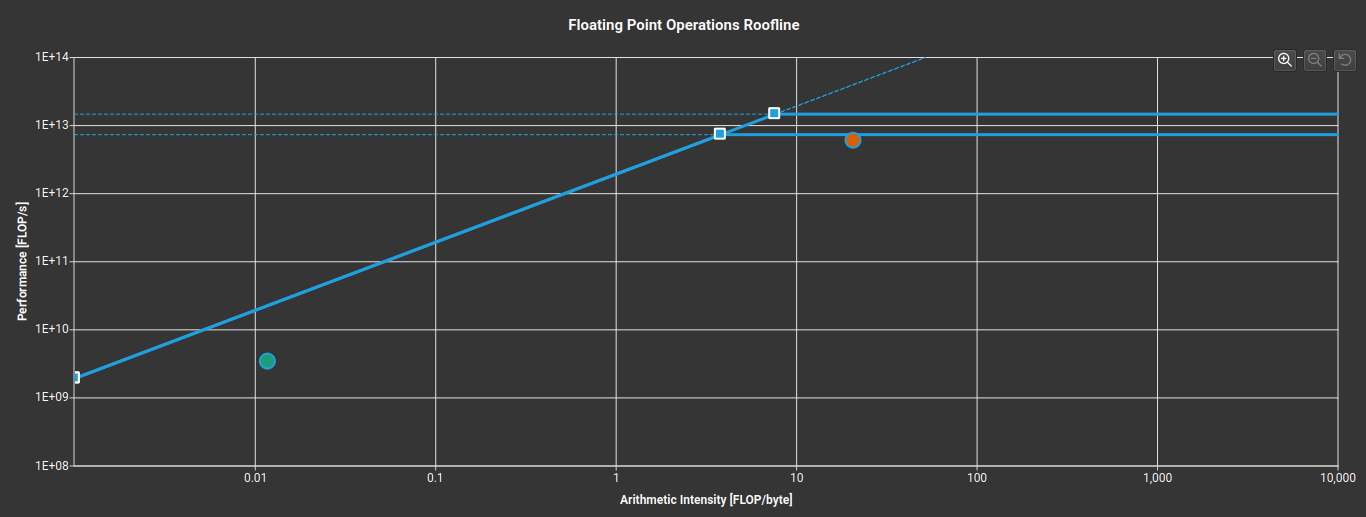}
    \caption{
        Solid Mechanics - NVIDIA A100 Roofline.
        \textmd{
            We observe that the NVIDIA A100 is computation bound, giving a potential explanation for the lower performance when compared to the AMD MI250X.
        }
    }
    \label{fig:roofline}
\end{figure}

In the case of the sole finite element operator, the AMD MI250X may outperform the NVIDIA A100 (see Figure~\ref{fig:SolidMechKernelThroughput}), showcasing the potential of the AMD MI250X for handling this specific arithmetically intense computation when the NVIDIA A100 is computation bound (see Figure~\ref{fig:roofline} ).
However, when considering the performance of the whole CG iteration, the situation is reversed due to the dot product implementation (see Figure~\ref{fig:SolidMechThroughput}) and the current superior ability of NVIDIA hardware to handle reduction algorithms.

This inversion of the higher-performing GPU between the sole operator and the whole CG iteration highlights the critical role of the dot product operation in determining the performance of the CG method. 
The observed discrepancy underscores the importance of optimizing the dot product implementation for both GPUs to ensure efficient performance and strong scalability in numerical simulations employing the CG method.

In our specific applications, we do not use the CG method as our primary solver, but rather employ other solvers tailored to the needs of the problems at hand. Given this context, we have chosen not to focus further on optimizing the CG method itself, such as by fusing operations into a single kernel. Instead, our attention is directed towards understanding and improving the performance of the individual components, such as the finite element operator, which are essential building blocks for various solvers.

\section{Conclusions}
\label{sec:conclusions}
In our study, we have developed a low-order matrix-free finite element implementation on GPUs that demonstrates competitive performance when compared to high-order methods on the same architecture.
While it is essential to acknowledge the inherent differences between low-order and high-order finite element methods, our implementation has successfully addressed many of the challenges associated with low-order methods on GPUs.
As a result, we have achieved performance levels that are comparable to those of high-order methods, highlighting the potential of our low-order approach for practical applications on GPU-based systems.
In conclusion, our study has demonstrated several key findings and highlighted potential future research directions:
\begin{itemize}
    \item Although the current work focuses on hexahedral elements, the described algorithm is well suited for simplicial elements.
    Investigating the performance and adaptability of the algorithm for simplices could further expand its applicability to a wider range of problems.

    \item Assessing the feasibility and performance of the matrix-free approach for quadratic and cubic finite elements would provide valuable insights into the benefits and limitations of the method for higher-order elements.

    \item Developing matrix-free preconditioners for matrix-free operators represents an exciting research direction that could significantly enhance the overall performance of the solvers employed in our applications.
\end{itemize}

\section{Acknowledgments}
\label{sec:acknowledgments}
This research was supported by the Exascale Computing Project (17-SC-20-SC), a collaborative effort of the U.S. Department of Energy Office of Science and the National Nuclear Security Administration.
Portions of this work were performed under the auspices of the U.S. Department of Energy by Lawrence Livermore National Laboratory under Contract DE-AC52-07-NA27344 (LLNL-CONF-847789).


\end{document}